\documentclass[12pt]{article}
 
\usepackage{amssymb,amsmath}
 
\usepackage{accents} 

\DeclareMathAlphabet{\eufrak}{U}{}{}{} 
\SetMathAlphabet\eufrak{normal}{U}{euf}{m}{n}
\SetMathAlphabet\eufrak{bold}{U}{euf}{b}{n}

\numberwithin{equation}{section}

\newenvironment{Proof}{\removelastskip\par\medskip
\noindent{\em Proof.} \rm}{\penalty-20\null\hfill$\square$\par\medbreak}

\allowdisplaybreaks

 \def\real{{\mathord{\mathbb R}}}
 
 \def\inte{{\mathord{\mathbb N}}}
 
 \def\qu{{\mathord{\mathbb Z}}}
 
\def\ind{{\bf 1}}

 \def\real{{\mathord{{\rm I\kern-3pt R}}}}        % Fake blackboard bold R.
 
 \def\inte{{\mathord{{\rm I\kern-3pt N}}}}
 \def\sZZ{{\rm Z\kern-.45em{}Z}}
 \def\sQQ{{\kern 0.27em \vrule height1.45ex width0.03em depth0em
           \kern-0.30em \rm Q}}
 \def\qu{{\mathchoice
         {\sQQ}
         {\sQQ}
   {\kern 0.225em \vrule height1.05ex width0.025em depth0em \kern-0.25em \rm Q}
   {\kern 0.180em \vrule height0.78ex width0.020em depth0em \kern-0.20em \rm Q}
         }}
 \def\sGG{{\kern 0.27em \vrule height1.45ex width0.03em depth0em
           \kern-0.30em \rm G}}
 \def\gg{{\mathchoice
         {\sGG}
         {\sGG}
   {\kern 0.225em \vrule height1.05ex width0.025em depth0em \kern-0.25em \rm G}
   {\kern 0.180em \vrule height0.78ex width0.020em depth0em \kern-0.20em \rm G}
         }}

 \newtheorem{prop}{Proposition}[section]
 \newtheorem{lemma}[prop]{Lemma}

 \newtheorem{theorem}[prop]{Theorem}

\def\E{\mathop{\hbox{\rm I\kern-0.20em E}}\nolimits}

 \newcounter{hyp}
\title{\huge Factorial moments of point processes} % Main lemmas} 
 
\author{
\large 
Jean-Christophe Breton\\ 
\small 
IRMAR - UMR CNRS 6625\\
\small
Universit\'e de Rennes 1\\
\small
Campus de Beaulieu\\
\small
F-35042 Rennes Cedex\\
\small
France\\
\and
 Nicolas Privault 
\\ 
\small
Division of Mathematical Sciences 
\\ 
\small 
School of Physical and Mathematical Sciences 
\\ 
\small
Nanyang Technological University 
\\ 
\small 
21 Nanyang Link 
\\ 
\small
Singapore 637371
}

\begin{document}

\hyphenation{func-tio-nals} 
\hyphenation{Privault} 

\maketitle 

\vspace{-0.9cm}

\baselineskip0.6cm
 
\begin{abstract} 
 We derive joint factorial moment identities 
 for point processes with Papangelou 
 intensities. 
 Our proof simplifies previous approaches 
 to related moment identities and includes 
 the setting of Poisson point processes. 
 Applications are given to random transformations 
 of point processes and to their distribution 
 invariance properties. 
\end{abstract} 
 
\noindent {\bf Key words:} 
 Point processes, 
 Papangelou intensities, 
 factorial moments, 
 moment identities. 
\\ 
{\em Mathematics Subject Classification (2010):} 60G57; 60G55; 60H07. 
 
\baselineskip0.7cm
 
\bigskip 
 
\section{Introduction} 
 Consider the compound Poisson random variable 
\begin{equation} 
\label{1} 
 \beta_1 Z_{\alpha_1} + \cdots + \beta_p Z_{\alpha_p} 
\end{equation} 
 where $\beta_1,\ldots ,\beta_p \in \real$ are constant 
 parameters and $Z_{\alpha_1}, \ldots , Z_{\alpha_p}$ is 
 a sequence of independent Poisson random variables 
 with respective parameters $\alpha_1,\ldots , \alpha_p \in \real_+$. 
 The L\'evy-Khintchine formula 
$$ 
 \Bbb{E} [ e^{t ( \beta_1 Z_{\alpha_1} + \cdots + \beta_p Z_{\alpha_p} )}] 
 = e^{ \alpha_1 ( e^{ \beta_1 t} -1 ) + \cdots 
  + \alpha_p ( e^{ \beta_p t} -1 )} 
$$ 
 shows that the cumulant of order $k\geq 1$ of \eqref{1} 
 is given by 
$$ 
 \alpha_1 \beta_1^k + \cdots + \alpha_p \beta_p^k. 
$$ 
 As a consequence of the Fa\`a di Bruno formula the 
 moment of order $n \geq 1$ of \eqref{1} is given by 
\begin{equation} 
\label{eq:2} 
 \Bbb{E}\left[ 
 \left( 
 \sum_{i=1}^p  \beta_i Z_{\alpha_i} 
 \right)^n 
 \right] 
 = 
 \sum_{m=1}^n 
 \sum_{ ~ P_1 \cup \cdots \cup P_m = \{ 1 , \ldots , n \} ~ } 
 \sum_{ i_1,\ldots ,i_m = 1}^p 
 \beta_{i_1}^{|P_1|} 
 \alpha_{i_1} 
 \cdots 
 \beta_{i_m}^{|P_m|} 
 \alpha_{i_m}, 
\end{equation} 
 where the above sum runs over all partitions 
 $P_1,\ldots ,P_m$ of $\{1,\ldots ,n\}$. 
\\ 
 
 Such cumulant-type moment identities have been 
 extended to Poisson stochastic integrals of 
 random integrands in 
 \cite{momentpoi} through the use of the Skorohod 
 integral on the Poisson space, cf. \cite{priinvcr}, \cite{prinv}. 
 The construction of the Skorohod integral has been extended to 
 point processes with 
 Papangelou intensities in \cite{torrisi2}, 
 and in \cite{flint}, the moment identities 
 of \cite{momentpoi} 
 have been extended to point processes with 
 Papangelou intensities via simpler 
 proofs based on an induction argument. 
\\ 
 
 In this paper we consider the factorial moments 
 of point processes and show that similar moment 
 identities can be deduced with an even simpler proof. 
 We apply those identities to random transformations 
 of Poisson processes and 
 point processes with Papangelou intensities, 
 and to their distribution invariance properties. 
\\ 
 
 Let $X$ be a Polish space equipped 
 with a $\sigma$-finite measure $\sigma(dx)$. 
 Let $\Omega^X$ denote 
 the space of configurations whose elements 
 $\omega \in \Omega^X$ are identified with the  Radon point measures 
 $\displaystyle \omega = \sum_{x\in \omega} \epsilon_x$, 
 where $\epsilon_x$ denotes the Dirac measure at $x\in X$. 
 A point process is a probability measure $P$ on $\Omega^X$ equipped 
 with the $\sigma$-algebra ${\cal F}$ generated by the topology of 
 vague convergence. 
 In the sequel for a (possibly random) set $A$ 
 we let 
 $N(A)(\omega ) =\int_X \ind_A(x)\ \omega(dx)$ 
 denote the cardinality of $\omega\cap A(\omega)$. 
\\
 
 Typically, for a Poisson point process with intensity $\sigma$, 
 $\omega(A)$ is distributed according to a Poisson distribution 
 with parameter $\sigma(A)$ for all (non random) $A\in{\cal F}$ 
 and $N (A)$ is independent of 
 $N (B)$ whenever $A, B\in{\cal F}$ are disjoint non random. 
\\ 
 
 Point processes can be 
 characterized by their Campbell measure $C$ defined on 
 ${\cal B}(X) \otimes {\cal F}$ by 
$$
C(A\times B) : =\Bbb{E}\left[\int_X \ind_A(x)\ind_B(\omega\setminus\{x\})\
\omega(dx)\right], 
 \qquad 
 A\in{\cal B}(X), \quad B\in{\cal F}. 
$$ 
 Recall the Georgii-Nguyen-Zessin identity 
\begin{equation}
\label{eq:GNZ0}
 \Bbb{E} \left[ 
 \int_X u ( x , \omega ) 
 \omega ( dx ) 
 \right] 
 = 
 \int_{\Omega^X} 
 \int_X u ( x, \omega \cup x ) 
C(dx,d\omega ) 
, 
\end{equation}
 for all measurable function $u: \Omega^X\times X\to \real$ such that both sides of \eqref{eq:GNZ0} make sense. 
 In particular, a Poisson point process with intensity $\sigma$ 
 is a point process with Campbell measure $C=\sigma\otimes P$, 
 and the Poisson measure with intensity $\sigma (dx)$ will be 
 denoted by $\pi_\sigma$. 
\\ 
 
 In the sequel, we consider Papangelou point processes, 
 i.e. point processes whose Campbell measure $C(dx,d \omega )$ is 
 absolutely continuous with respect to 
 $\sigma\otimes P$, i.e. 
$$ 
 C(dx,d\omega ) = c(x,\omega ) \sigma (dx ) P(d\omega ), 
$$ 
 where the density $c(x,\omega)$ is called the Papangelou density. 
 In this case the identity \eqref{eq:GNZ0} reads 
\begin{equation}
\label{eq:GNZ}
 \Bbb{E} \left[ 
 \int_X u ( x, \omega ) 
 \omega ( dx ) 
 \right] 
 = 
 \Bbb{E} \left[ 
 \int_X u ( x , \omega \cup x ) 
 c(x,\omega ) 
 \sigma ( dx ) 
 \right] 
, 
\end{equation}
and $c(x, \omega)=1$ for Poisson point process with intensity $\sigma$. 
\\
 
 This paper is organized as follows.
 In Section \ref{sec:fac-moment}, we derive factorial moment identities 
 for random point measure of random sets 
 in Propositions~\ref{dldd1} and \ref{klddd2}, 
 and in Section~\ref{3.1} we apply those identities 
 to point process transformations in Proposition~\ref{djhjklddd111}. 
 In Section \ref{sec:moment}, we show that the corresponding 
 moment identities can be recovered by combinatorial arguments, 
 cf. Proposition~\ref{djld}. 
 In Section \ref{sec:invariance}, we recover some recent 
 results on the invariance of Poisson random measures 
 under interacting transformations, with simplified proofs.  
\section{Factorial moments} 
\label{sec:fac-moment}
 Let 
$$ 
 x_{(n)}=x(x-1) \cdots (x-n+1), \qquad 
 x\in \real, \quad n \in \inte, 
$$ 
 denote the falling factorial product. 
We are interested in the {\it factorial moments} 
 $\mu_n^f(N(A))=\Bbb{E}[N(A)_{(n)}]$ when $N(A)$ 
 is the random point measure of a random set $A$. 
 Denoting by $\Omega_0^X$ the set of finite configurations 
 in $\Omega^X$, the compound Campbell density 
$$ 
\hat{c}: \Omega_0^X\times\Omega^X \longrightarrow \real_+ 
$$ 
 is defined inductively by
\begin{equation}
\label{eq:hatc}
 \hat{c} ( 
 \{x_1,\ldots ,x_n , y \} , 
 \omega ) 
 : = 
 c( y , \omega ) 
 \hat{c} ( 
 \{x_1,\ldots ,x_n \} , 
 \omega \cup \{ y\} ) 
, \qquad n \geq 0. 
\end{equation} 
 Given $\eufrak{x}_n =(x_1, \ldots, x_n)\in X^n$, we will use 
 the notation $\varepsilon^+_{\eufrak{x}_n }$ 
 for the operator 
$$ 
 ( \varepsilon^+_{\eufrak{x}_n } F)(\omega)=F(\omega\cup\{x_1, \ldots, x_n\}), 
 \qquad 
 \omega \in \Omega, 
$$ 
 where $F$ is any random variable on $\Omega^X$. 
 With this notation we also have 
$$ 
\hat{c} ( \eufrak{x}_n ,\omega )
 = 
 c(x_1, \omega) c(x_2, \omega\cup\{x_1\})
c(x_3, \omega \cup\{x_1, x_2\}) \cdots
c(x_n, \omega \cup\{x_1, \ldots, x_{n-1}\})
. 
$$ 
 In addition, we define the random measure 
 $\hat{\sigma}^n ( d\eufrak{x}_n )$ on $X^n$ by 
$$ 
 \hat{\sigma}^n ( d\eufrak{x}_n ) 
 = 
 \hat{c} ( \eufrak{x}_n ,\omega ) 
 \sigma^n ( d\eufrak{x}_n ) 
 = 
 \hat{c} ( \eufrak{x}_n ,\omega ) 
 \sigma (dx_1) \cdots \sigma (dx_n) 
, 
$$ 
 with 
$ 
 \sigma^n ( d\eufrak{x}_n ) 
 = 
 \sigma (dx_1) \cdots \sigma (dx_n) 
$. 
\begin{prop} 
\label{klddd2} 
 Let $A=A (\omega )$ be a random set.
 For all $n\geq 1$ and sufficiently integrable random variable $F$, 
 we have 
\begin{equation*} 
\label{djkldd10} 
\Bbb{E} \left[ 
F\  N(A)_{(n)} 
\right] 
= 
 \Bbb{E} \left[ 
 \int_{X^n} 
 \varepsilon^+_{\eufrak{x}_n}( 
 F \ind_{A^n} ( x_1, \ldots , x_n )  )(\omega)
\ 
 \hat{\sigma}^n ( dx_1, \ldots , dx_n ) 
 \right] 
. 
\end{equation*} 
\end{prop} 
\begin{Proof} 
 We show by induction on $n\geq 1$ that
\begin{equation} 
\label{eq:step1}
\Bbb{E} \left[ 
F\  N(A)_{(n)} 
\right] 
= 
 \Bbb{E} \left[ 
 \int_{X^n} 
 \varepsilon^+_{\eufrak{x}_n}( 
 F  \ind_{A} (x_1) \cdots  \ind_{A} (x_n)  )(\omega)
\ 
 \hat{\sigma}^n ( d\eufrak{x}_n ) 
 \right]. 
\end{equation} 
 Clearly the formula 
$$
 \Bbb{E} \left[ 
 F N ( A ) 
 \right] 
 = 
 \Bbb{E} \left[ 
 \int_{X} 
 \varepsilon^+_{x}( 
 F  \ind_{A} (x) )(\omega)\ 
c( x , \omega )\ 
 \sigma (dx) \right] 
$$ 
 holds at the rank $n=1$ due to \eqref{eq:GNZ} 
 applied to $ u( x , \omega )=F(\omega)\ind_{A(\omega)}(x)$. 
 Next, assuming that \eqref{eq:step1} holds at the rank $n$, 
 we apply it with $F$ replaced by $F(N(A)-n)$ and get 
\begin{eqnarray} 
\nonumber 
\lefteqn{ 
 \! \! \! \! \! \! \! \! \! \! \! \! \! \! \! \! \! \! \! \! \! \! \! \! \! \! 
 \Bbb{E} \big[ 
 F\ N(A) \cdots (N(A)- n ) 
 \big] 
 = 
 \Bbb{E} \left[ 
 \int_{X^n} 
\varepsilon^+_{\eufrak{x}_n} ( 
 F\ (N(A)-n ) 
 \ind_{A} (x_1) 
 \cdots 
 \ind_{A} (x_n) 
 )(\omega) \
 \hat{\sigma}^n ( d\eufrak{x}_n ) 
 \right] 
} 
\\ 
\nonumber 
 & = & 
 \Bbb{E} \left[ 
 \int_{X^n} 
 \varepsilon^+_{\eufrak{x}_n}( 
 F\ N(A) 
 \ind_{A} (x_1) 
 \cdots 
 \ind_{A} (x_n) 
 )  (\omega) \ 
 \hat{\sigma}^n ( d\eufrak{x}_n ) 
 \right] 
\\ 
\nonumber 
 & & 
 - n 
 \Bbb{E} \left[ 
 \int_{X^n} 
\varepsilon^+_{\eufrak{x}_n} ( 
 F 
 \ind_{A} (x_1) 
 \cdots 
 \ind_{A} (x_n) 
 ) 
 (\omega)  \
 \hat{\sigma}^n ( d\eufrak{x}_n ) 
 \right] 
\\ 
\label{eq:step1b}
 & = & 
 \Bbb{E} \left[ 
 \int_{X^n} 
N(\varepsilon^+_{\eufrak{x}_n}(A))(\omega) 
 \varepsilon^+_{\eufrak{x}_n}( 
 F \ind_{A} (x_1) 
 \cdots 
 \ind_{A} (x_n) 
 ) 
 (\omega) \
 \hat{\sigma}^n ( d\eufrak{x}_n ) 
 \right] 
, 
\end{eqnarray}
 where in \eqref{eq:step1b} we used the relation 
\begin{eqnarray*}
\varepsilon^+_{\eufrak{x}_n}(N(A))(\omega)
&=&N(\varepsilon^+_{\eufrak{x}_n}(A))(\omega)+\sum_{i=1}^n \delta_{x_i}(\varepsilon^+_{\eufrak{x}_n}(A)(\omega))\\
&=&N(\varepsilon^+_{\eufrak{x}_n}(A))(\omega)+\sum_{i=1}^n \varepsilon^+_{\eufrak{x}_n}(\ind_{A}(x_i))(\omega). 
\end{eqnarray*}
 Next, with $\eufrak{x}_{n+1}=(x_1,\ldots , x_n,x_{n+1})$, 
recalling that $N(A)=\int_X \ind_A(x)\ \omega(dx)$ 
 and applying \eqref{eq:GNZ} to 
\begin{eqnarray*} 
 u( x , \omega ) &: =& 
 \varepsilon^+_{\eufrak{x}_n}( 
 F \ind_{A} (x_1) 
 \cdots 
 \ind_{A} (x_n) \ind_A(x)
 ) 
 (\omega)  
\ \hat c({\eufrak{x}_n},\omega)
\\ 
&=& 
 F(\omega\cup\{x_1,\ldots, x_n\}) \ind_{A(\omega\cup\{x_1,\ldots, x_n\} )} (x_1) 
 \cdots \ind_{A(\omega\cup\{x_1,\ldots, x_n\})} (x_n)\\
&& \hskip 6cm \times\ind_{A(\omega\cup\{x_1,\ldots, x_n\})}(x)
\ \hat c({\eufrak{x}_n},\omega)
\end{eqnarray*} 
 for fixed $x_1, \ldots, x_n$ 
 with the relation 
 $\varepsilon^+_{\eufrak{x}_{n+1}}=\varepsilon^+_{x_{n+1}}\circ\varepsilon^+_{\eufrak{x}_n}$ 
 we find 
\begin{eqnarray*}
\nonumber 
\lefteqn{ 
\Bbb{E} \left[ 
 F\ N(A)\cdots (N(A)- n ) 
 \right] 
} 
\\
&=& 
 \label{eq:tech23}
 \Bbb{E} \left[ 
 \int_{X^{n+1}} 
 \varepsilon^+_{\eufrak{x}_{n+1}}\big( 
 F 
 \ind_{A} (x_1) 
 \cdots 
 \ind_{A} (x_{n+1} ) \big)  
 (\omega) 
\ \hat{c} (\eufrak{x}_n ,\omega \cup \{ x_{n+1}\} )
\ c(x_{n+1},\omega)  
\right. 
\\
\nonumber
 &&\hskip 7cm \sigma (dx_1) \cdots  \sigma (dx_{n+1} ) \Big] 
\\ 
\nonumber 
 & = & 
 \Bbb{E} \left[ 
 \int_{X^{n+1}} 
\varepsilon^+_{\eufrak{x}_{n+1}} ( 
 F 
 \ind_{A} (x_1) 
 \cdots 
 \ind_{A} (x_{n+1} ) ) 
 ( \omega) 
\ \hat{\sigma}^n ( d\eufrak{x}_{n+1} ) 
 \right], 
\end{eqnarray*} 
 where on the last line we used \eqref{eq:hatc}.
\end{Proof} 
\noindent 
 By induction, 
 in the next Proposition~\ref{dldd1} we also obtain 
 a joint factorial moment identity for a.s. 
 disjoint (random) sets $A_1, \ldots, A_p$. 
 It extends the classical identity 
\begin{equation} 
\label{djklddda} 
 \Bbb{E} \left[ 
 N(A_1)_{(n_1)}  \cdots 
N(A_p)_{(n_p)}  \right] 
 = 
 \int_{A_1^{n_1} \times \cdots 
 \times A_p^{n_p} } 
 \rho_n (x_1,\ldots ,x_n)\ 
 \sigma (dx_1) \cdots  \sigma (dx_n) 
, 
\end{equation}  
 for deterministic disjoint sets $A_1, \ldots , A_p$, 
 where  $\rho_n (x_1,\ldots ,x_n)$ 
 is the correlation function of 
 the point process and $n=n_1+\cdots+n_p$.  
\begin{prop} 
\label{dldd1} 
 Let $n = n_1 + \cdots + n_p$ and $A_1(\omega), \ldots, A_p(\omega)$ 
 be measurable and disjoint for almost all $\omega\in\Omega$, then 
\begin{equation} 
\label{djkldd11} 
 \Bbb{E} \left[ 
 F\ N(A_1)_{(n_1)} 
 \cdots N(A_p)_{(n_p)} 
 \right] 
 = 
 \Bbb{E} \left[ 
 \int_{X^n} 
 \varepsilon^+_{\eufrak{x}_n}( 
 F 
(\ind_{A_1^{n_1}}\otimes \cdots \otimes
 \ind_{A_p^{n_p}})(\eufrak{x}_n ) 
 ) \ \hat{\sigma}^n ( d\eufrak{x}_n ) 
 \right] 
. 
\end{equation} 
\end{prop}
\begin{Proof} 
 We proceed by induction on $p \geq 1$. 
 For $p=1$, the identity reduces to that of Proposition \ref{klddd2}. 
 We assume that the identity holds true for $p$ and show it for $p+1$. 
 Let $n=n_1+\cdots+n_p$ and $m=n+n_{p+1}$, we have: 
\begin{eqnarray} 
\nonumber 
\lefteqn{ 
 \Bbb{E} \left[ 
 F\ N(A_1)_{(n_1)} 
 \cdots N(A_{p+1})_{(n_{p+1})} 
 \right] 
} 
\\ 
\nonumber 
& = & 
 \Bbb{E} \left[ 
 \int_{X^n} 
 \varepsilon^+_{\eufrak{x}_n} 
 \left( 
 F N(A_{p+1})_{(n_{p+1})}
 (\ind_{A_1^{n_1}} \otimes \cdots \otimes
 \ind_{A_p^{n_p}}) (x_1,\ldots, x_n) 
 \right) 
 \hat{\sigma}^n ( d\eufrak{x}_n ) 
 \right] 
\\
 \label{eq:tech31} 
& = & 
 \Bbb{E} \left[ 
 \int_{X^n} 
 N(\varepsilon^+_{\eufrak{x}_n} 
 (A_{p+1})_{(n_{p+1})})
 \varepsilon^+_{\eufrak{x}_n}( 
 F 
 (\ind_{A_1^{n_1}}\otimes \cdots\otimes \ind_{A_p^{n_p}})
(x_1, \ldots, x_n) 
 ) 
\ \hat{\sigma}^n ( d\eufrak{x}_n ) 
 \right] 
\\
\nonumber
 & = & 
 \Bbb{E} \left[ 
 \int_{X^n}\int_{X^{n_{p+1}}} 
 \varepsilon^+_{\eufrak{y}_{n_{p+1}}} 
 \Big( 
 \varepsilon^+_{\eufrak{x}_n} 
 \left( 
 F 
 \ind_{A_1^{n_1}} \otimes\cdots\otimes\ind_{A_p^{n_p}} 
 (x_1, \ldots,x_n) 
 \right) \hat{c} ( \{x_1,\ldots ,x_n \} , \omega )
\right . 
\\
\nonumber 
&& \ind_{\varepsilon^+_{\eufrak{x}_n}(A_{p+1}^{n_{p+1}})}(y_1, \ldots, y_{n_{p+1}})\Big)
\ \hat{c} ( \{y_1,\ldots ,y_{n_{p+1}} \} , \omega )\
 \sigma (dy_1) \cdots  \sigma (dy_{n_{p+1}})
\sigma (dx_1) \cdots  \sigma (dx_n) \Big] 
\\
\label{eq:tech32}
&&
\\
\label{eq:tech33}
 & = & 
 \Bbb{E} \left[ 
 \int_{X^m} 
 \varepsilon^+_{\eufrak{x}_m} \left( 
 F 
 \ind_{A_1^{n_1}}\otimes \cdots\otimes \ind_{A_{p+1}^{n_{p+1}}} 
 (x_1,\ldots, x_m) \right)
 \hat{\sigma} ( dx_1,\ldots , dx_m ) 
 \right],   
\end{eqnarray} 
 where in \eqref{eq:tech31} we used
\begin{eqnarray*}
\varepsilon^+_{\eufrak{x}_n}(N(A_{p+1})_{(n_{p+1})}) 
&=&\left(N(\varepsilon^+_{\eufrak{x}_n}(A_{p+1}))+\sum_{i=1}^n \delta_{x_i}(\varepsilon^+_{\eufrak{x}_n}(A_{p+1}))\right)_{(n_{p+1})} 
\\
&=&\left(N(\varepsilon^+_{\eufrak{x}_n}(A_{p+1}))+\sum_{i=1}^n \varepsilon^+_{\eufrak{x}_n}(\ind_{A_{p+1}}(x_i))\right)_{(n_{p+1})} 
\end{eqnarray*} 
 and observe that the contribution of the sum is zero since, for all $1\leq k\leq p$ and $1\leq i\leq n$, $\ind_{A_{p+1}}(x_i)\ind_{A_k}(x_i)=0$. 
 In \eqref{eq:tech32}, we noted $\eufrak{y}_{n_{p+1}}=(y_1, \ldots, y_{n_{p+1}})$ and used Proposition~\ref{klddd2} with, for a fixed $\eufrak{x}_n=(x_1, \ldots, x_n)$,
$$
\widetilde F(\omega)=\varepsilon^+_{\eufrak{x}_n}\big( 
 F (\ind_{A_1^{n_1}}\otimes \cdots \otimes\ind_{A_p^{n_p}})
(x_1,\ldots,x_n) \big)(\omega) \ \hat{c}(\{x_1,\ldots,x_n\},\omega)) 
$$
and the set $\varepsilon^+_{\eufrak{x}_n}(A_{p+1})$. 
Finally in \eqref{eq:tech33}, we used the following consequence of \eqref{eq:hatc}
$$
\hat{c} \big( \{x_1,\ldots ,x_n \} , \omega\cup\{y_1, \ldots, y_{n_p}\}\big)
\ \hat{c} \big( \{y_1,\ldots ,y_{n_p} \} , \omega\big)
=\hat{c}\big( \{x_1,\ldots ,x_n, y_1, \ldots, y_{n_p}\} , \omega\big)
$$
together with $\varepsilon_{\eufrak{y}_{n_p}}^+\circ \varepsilon_{\eufrak{x}_{n}}^+
=\varepsilon_{\eufrak{y}_{n_p}\cup\eufrak{x}_n}^+$.
\end{Proof} 
\section{Transformations of point processes} 
\label{3.1} 
 Consider the finite difference operator
$$
D_xF(\omega)=F(\omega\cup\{x\})-F(\omega)
$$ 
 where $F$ is any random variable on $\Omega^X$. 
Note that multiple finite difference operator expresses
\begin{equation}
\label{eq:Dmult}
D_\Theta F=\sum_{\eta\subset\Theta} (-1)^{|\Theta|+1+|\eta|} F(\omega\cup\eta)
\end{equation}
where the summation above holds over all (possibly empty) 
 subset $\eta$ of $\Theta$. 
 Let ${\eufrak{x}_{n}}=\{x_1, \ldots, x_n\}$, from the relation 
\begin{eqnarray} 
\nonumber
 \varepsilon^+_{\eufrak{x}_n } 
 ( 
 u_1 ( x_1 , \omega )  
 \cdots 
 u_n ( x_n , \omega )  
 ) 
&=& 
 \varepsilon_{x_1,\ldots ,x_n}^+ 
 ( 
 u_1 ( x_1, \omega ) 
 \cdots 
 u_n ( x_n , \omega ) 
 ) 
\\
\label{eq:espilon-mult}
&=& 
 \sum_{\Theta \subset \{1,\ldots ,n\} } 
 D_{\Theta} 
 \big(  u_1 ( x_1, \omega ) \cdots  u_n ( x_n , \omega )\big) 
, 
\end{eqnarray}
 where 
 $D_{\Theta} = D_{x_1} \cdots D_{x_l}$ when 
 $\Theta = \{ 1, \ldots , l\}$
and from \eqref{djkldd11} we have 
\begin{eqnarray} 
\nonumber
\lefteqn{ 
 \! \! \! \! \! \! \! \! \! \! \! \! \! \! \! \! \! \! \! \! \! \! \! \! \! \! \! 
 \Bbb{E} \left[ 
 F\ N(A_1)_{(n_1)} 
 \cdots N(A_p)_{(n_p)} 
 \right] 
 = 
 \Bbb{E} \left[ 
 \int_{X^n} 
 \varepsilon^+_{\eufrak{x}_n}( 
 F 
 (\ind_{A_1^{n_1}}\otimes \cdots \otimes
 \ind_{A_p^{n_p}})(\eufrak{x}_n ) 
 ) 
\ {\hat\sigma}^n ( d\eufrak{x}_n ) 
 \right] 
} 
\\
\label{eq:DThetajoint}
&=&\sum_{\Theta \subset \{1,\ldots ,n\} } 
\Bbb{E} \left[ 
 \int_{X^n} 
 D_{\Theta}( 
 F 
(\ind_{A_1^{n_1}}\otimes \cdots \otimes
 \ind_{A_p^{n_p}})(\eufrak{x}_n ) 
 ) 
\ {\hat\sigma}^n ( d\eufrak{x}_n ) 
 \right] 
\end{eqnarray} 
 for $n_1+\cdots +n_p=n$ and a.s. disjoint sets 
 $A_1(\omega), \ldots, A_p(\omega)$. 
 The next lemma will be useful in Proposition~\ref{djhjklddd111} 
 to characterize 
 the invariance of transformations of point processes 
 from \eqref{eq:DThetajoint}.
\begin{lemma} 
\label{l12.32} 
 Let $m\geq 1$ and assume that 
 for all $x_1,\ldots ,x_m \in X$ 
 the processes $u_i : \Omega^X \times X \longrightarrow \real$, $1\leq i\leq m$ satisfy the condition 
\begin{equation} 
\label{djdkd1} 
 D_{\Theta_1} 
 u_1 ( x_1 , \omega ) 
 \cdots 
 D_{\Theta_m} 
 u_m ( x_m , \omega ) 
 = 
 0, 
\end{equation} 
 for every family $\{ \Theta_1, \ldots , \Theta_m \}$ 
 of (non empty) subsets such that 
 $\Theta_1 \cup \cdots \cup \Theta_m = \{1,\ldots ,m \}$, 
 for all $x_1,\ldots ,x_m \in X$ and all $\omega\in\Omega^X$. 
 Then we have 
\begin{equation}
\label{eq:DDD}
 D_{x_1} \cdots D_{x_m}\big( u_1(x_1,\omega)  \cdots u_m(x_m,\omega)\big)=0 
\end{equation}
 for all $x_1,\ldots ,x_m \in X$ and all $\omega\in\Omega^X$. 
\end{lemma} 
\begin{Proof} 
 It suffices to note that 
\begin{equation}
\label{eq:DDD0}
 D_{x_1} \cdots D_{x_n} \big( u_1( x_1 , \omega ) 
 \cdots u_l ( x_l , \omega )\big) =
 \sum_{\Theta_1\cup\cdots\cup \Theta_l = \{1, \ldots, n\}} D_{\Theta_1} 
 u_1( x_1 , \omega ) \cdots D_{\Theta_l} u_l( x_l , \omega ),
\end{equation} 
 where the above sum is not restricted to partitions, 
 but includes all (possibly empty) sets 
 $\Theta_1, \ldots , \Theta_l$ whose union is 
 $\{1, \ldots, n\}$. 
\end{Proof} 
 In the next result, we recover Theorem~5.1 of \cite{flint} in  
 a more direct way due to the use of factorial moments, 
 but using a different cyclic type condition. 
 Condition~\eqref{djdkd1.0} below is interpreted 
 by saying that 
\begin{equation} 
\label{eq:interpretation} 
 D_{\Theta_1} 
 h_1 ( \tau ( x_1 , \omega ) ) 
 \cdots 
 D_{\Theta_m} 
 h_m ( \tau ( x_m , \omega ) ) 
 = 
 0, 
\end{equation} 
 for any family $h_1,\ldots ,h_m$ of bounded 
 real-valued Borel functions on $Y$. 
\begin{prop} 
\label{djhjklddd111} 
 Let $\tau : \Omega^X \times X \to Y$ 
 be a random transformation such that $\tau ( \cdot , \omega ): X \to Y$ 
 maps bijectively $\sigma$ to $\mu$ for all $\omega \in \Omega^X$, 
 i.e. 
$$ 
\sigma\circ\tau( \cdot , \omega )^{-1} = \mu, 
 \qquad 
 \omega \in \Omega^X, 
$$ 
 and satisfying the condition 
\begin{equation} 
\label{djdkd1.0} 
 D_{\Theta_1} 
 \tau ( x_1 , \omega ) 
 \cdots 
 D_{\Theta_m} 
 \tau ( x_m , \omega ) 
 = 
 0, 
\end{equation} 
 for every family $\{ \Theta_1, \ldots , \Theta_m \}$ 
 of (non empty) subsets such that 
 $\Theta_1 \cup \cdots \cup \Theta_m = \{1,\ldots ,m \}$, 
 for all $x_1,\ldots ,x_m \in X$ and all $\omega\in\Omega^X$, 
 $m\geq 1$. 
 Then $\tau_* : \Omega^X \to \Omega^Y$ defined by 
\begin{equation*}
\label{eq:tau*}
 \tau_* \omega = \sum_{x\in \omega } \epsilon_{\tau ( x , \omega )}=\omega\circ \tau(\cdot,\omega)^{-1}, 
 \qquad 
 \omega \in \Omega^X, 
\end{equation*} 
 transforms a point process $\xi$ with Papangelou intensity 
 $c ( x , \omega )$ with respect to $\sigma\otimes P$ into a
 point process on $Y$ with correlation function 
$$ 
 \rho_\tau (y_1,\ldots ,y_n) 
 = 
 \Bbb{E} \big[ 
 \hat{c} ( 
 \{ \tau^{-1} ( y_1 , \omega ),\ldots ,
 \tau^{-1} ( y_n , \omega )\} , 
 \omega ) \big], 
$$ 
 $y_1,\ldots,y_n \in Y$, with respect to $\mu$. 
\end{prop} 
\begin{Proof} 
 Consider $B_1,\ldots, B_p$ disjoint {\em deterministic} 
 subsets of $Y$ such that $\mu(B_1), \ldots , \mu (B_p)$ 
 are finite. 
From interpretation \eqref{eq:interpretation}, 
Condition \eqref{djdkd1.0} ensures \eqref{djdkd1} for 
$u_k(x,\omega)=\ind_{B_{i_k}}(\tau(x,\omega))$ and, in turn, Lemma~\ref{l12.32} shows that 
\begin{equation}
\label{eq:cons1lemma}
 D_{x_1} \cdots D_{x_k} \big( 
 {\bf 1}_{B_{i_1}} ( \tau ( x_1 , \omega )) 
 \cdots 
 {\bf 1}_{B_{i_k}} ( \tau ( x_k , \omega )) 
 \big) = 0, 
 \qquad 
 x_1, \ldots , x_k \in X, 
\end{equation}
 for all $i_1,\ldots ,i_k \in \{ 1,\ldots ,p \}$ and 
 $\omega \in \Omega$. 
For $i=1, \ldots, p$, let $A_i(\omega)=\tau( \cdot , \omega )^{-1}(B_i)$ 
and $\tau_*N(B_i)$ be the cardinal of $\tau_*\omega\cap B_i$,  
 i.e. 
 $$
 \tau_*N(B_i)= \sum_{x\in\omega}\epsilon_{\tau(x , \omega )}(B_i)
 = \sum_{x\in\omega}\epsilon_{x}(A_i)
 = N(A_i).
$$ 
Then applying \eqref{eq:DThetajoint} with $F=1$ and the disjoint random sets 
$A_i(\omega)$, $i=1, \ldots, p$, yields
\begin{eqnarray} 
\nonumber
\lefteqn{
 \Bbb{E} \left[  \tau_*N(B_1)_{(n_1)} \cdots \tau_* N(B_p)_{(n_p)} \right] 
 =  \Bbb{E} \left[ N(A_1)_{(n_1)} \cdots N(A_p)_{(n_p)}  \right]
}\\
\nonumber &=& 
\sum_{\Theta \subset \{1,\ldots ,n\} } 
\Bbb{E} \left[ 
 \int_{X^n} 
 D_{\Theta} 
 ( 
 ( 
 \ind_{A_1^{n_1}}\otimes \cdots \otimes
 \ind_{A_p^{n_p}} ) (\eufrak{x}_n ) 
 )\ 
 {\hat\sigma}^n ( d\eufrak{x}_n ) 
 \right] 
\\
\nonumber
&=& 
\Bbb{E} \left[ 
 \int_{X^n} 
 D_{x_1} \cdots D_{x_n} 
 ( 
 ( 
 \ind_{B_1^{n_1}}\otimes \cdots \otimes
 \ind_{B_p^{n_p}})(\tau( \eufrak{x}_n , \omega ) ) 
 )\   
 {\hat\sigma}^n ( d\eufrak{x}_n ) 
 \right]
\\ 
\nonumber
 & & + 
\sum_{\Theta \subsetneq \{1,\ldots ,n\} } 
\Bbb{E} \left[ 
 \int_{X^n} 
 D_{\Theta} 
 ( 
 (\ind_{B_1^{n_1}}\otimes \cdots \otimes
 \ind_{B_p^{n_p}})(\tau( \eufrak{x}_n , \omega )) 
 )\ 
 {\hat\sigma}^n ( d\eufrak{x}_n ) 
 \right]
\\ 
\label{eq:generic0}
 & = & 
\sum_{\Theta \subsetneq \{1,\ldots ,n\} } 
\Bbb{E} \left[ 
 \int_{X^n} 
 D_{\Theta} 
 ( 
 (\ind_{B_1^{n_1}}\otimes \cdots \otimes
 \ind_{B_p^{n_p}})(\tau( \eufrak{x}_n , \omega )) 
 )\
 {\hat\sigma}^n ( d\eufrak{x}_n ) 
 \right]
\end{eqnarray} 
 $n_1+\cdots +n_p=n \geq 1$, 
 where $\tau(\eufrak{x}_n , \omega )$ stands for $\big(\tau(x_1 , \omega ), \ldots, \tau(x_n , \omega )\big)$ 
 and where we used \eqref{eq:cons1lemma}. 
 Next, without loss of generality the generic term of \eqref{eq:generic0} can be reduced to 
 the term with $\Theta = \{ 1 , \ldots , n-1 \}$ and using \eqref{eq:DDD0}, we have 
\begin{eqnarray*} 
\nonumber
\lefteqn{ 
 \Bbb{E} \left[ 
 \int_{X^n} 
 D_{x_1} \cdots D_{x_{n-1}} 
 ( 
 (\ind_{B_1^{n_1}}\otimes \cdots \otimes
 \ind_{B_p^{n_p}})(\tau(\eufrak{x}_n , \omega )) 
 ) 
 {\hat \sigma}^n ( d\eufrak{x}_n ) 
 \right] 
} 
\\ 
 & = & 
 \Bbb{E} \left[ 
 \int_X 
 \int_{X^{n-1}} 
 D_{x_1} \cdots D_{x_{n-1}} 
 ( 
 (\ind_{B_1^{n_1}}\otimes \cdots \otimes
 \ind_{B_p^{n_p}})\big(\tau( \eufrak{x}_{n-1}), y_n , \omega \big) 
\right. 
\\ 
 & & 
\hskip 2cm \hat{c} ( 
 \{ x_1 
 ,\ldots ,
 x_{n-1} 
 , 
 \tau^{-1} ( y_n , \omega ) 
 \} , 
 \omega ) ) \  
 {\sigma}^{n-1} ( d\eufrak{x}_{n-1} ) 
 \mu( dy_n) 
 \Big] 
\end{eqnarray*} 
 with the change of variable $y_n=\tau( x_n , \omega )$. 
 Finally, by applying the above argument recursively we obtain that 
\begin{eqnarray*} 
\lefteqn{ 
 \Bbb{E} \left[  \tau_*N(B_1)_{(n_1)} \cdots \tau_*N(B_p)_{(n_p)}  \right] 
} 
\\ 
 & = & 
 \int_{X^n} 
 (\ind_{B_1^{n_1}}\otimes \cdots \otimes
 \ind_{B_p^{n_p}})(\eufrak{y}_n) \
 \Bbb{E} \big[ 
 \hat{c} ( 
 \{ \tau^{-1} ( y_1 , \omega ) 
 ,\ldots ,
 \tau^{-1} ( y_n , \omega ) 
 \} , 
 \omega )\big]\ 
 {\mu}^n ( d\eufrak{y}_n ), 
\end{eqnarray*} 
 $n_1,\ldots ,n_p = n \geq 1$, 
 which recovers the definition of the correlation function of $\tau_*\xi$ (see \eqref{djklddda}). 
\end{Proof} 
 The proof of Proposition~\ref{djhjklddd111} also 
 shows that if 
 $A_1,\ldots, A_p$ are disjoint random 
 subsets of $X$ such that 
$$ 
 D_{\Theta_1} 
 {\bf 1}_{A_1(\omega )} ( x_1 ) 
 \cdots 
 D_{\Theta_m} 
 {\bf 1}_{A_m (\omega )} ( x_m ) 
 = 
 0, 
$$ 
 for every family $\{ \Theta_1, \ldots , \Theta_m \}$ 
 of (non empty) subsets such that 
 $\Theta_1 \cup \cdots \cup \Theta_m = \{1,\ldots ,m \}$, 
 for all $x_1,\ldots ,x_m \in X$ and all $\omega\in\Omega^X$, 
 $m \geq 1$, then we have 
$$ 
 \Bbb{E} \left[ N(A_1)_{(n_1)} \cdots N(A_p)_{(n_p)}  \right]
 = 
 \Bbb{E} \left[ 
 \int_{X^n} 
 ( 
 \ind_{A_1^{n_1}}\otimes \cdots \otimes
 \ind_{A_p^{n_p}} ) (\eufrak{x}_n ) 
 {\hat\sigma}^n ( d\eufrak{x}_n ) 
 \right] 
, 
$$ 
 $n_1+\cdots +n_p=n$. 
\subsubsection*{Example} 
 We consider an example of transformation 
 satisfying Condition~\eqref{djdkd1.0}, 
 based on conditioning by a random boundary, 
 more precisely conditioned by the random boundary of a convex Poisson hull. 
 We let $X = \real^d$ with norm $\Vert \cdot \Vert$ 
 and for all $\omega \in \Omega$ we denote by 
 $\omega_e \subset \omega$ the extremal vertices of 
 the convex hull of $\omega \cap B(0,1)$. 
 We also denote by ${\cal C} (\omega )$ the convex hull of $\omega$, 
 and denote $\accentset{\circ}{\cal C} (\omega )$ its interior. 
\\ 
 
 Consider a mapping ${\tau} : \Omega^X \times X \longrightarrow X$ 
 such that for all $\omega \in \Omega$, 
 ${\tau} ( \cdot , \omega ) : X \longrightarrow X$ leaves 
 $X \setminus \accentset{\circ}{\cal C} (\omega_e )$ invariant 
 (thus including the extremal vertices $\omega_e$ of ${\cal C} (\omega_e )$) 
 while the points inside $\accentset{\circ}{\cal C}(\omega_e )$ are shifted 
 depending on the data of $\omega_e$, i.e. we have 
\begin{equation} 
\label{adfdsfg} 
 {\tau} ( x , \omega ) = 
 \left\{ 
 \begin{array}{ll} 
 {\tau} ( x , \omega_e ), 
 & 
 x\in \accentset{\circ}{\cal C}(\omega_e ), 
\\ 
\\ 
 x , 
 & 
 x\in X\setminus \accentset{\circ}{\cal C}(\omega_e ). 
\end{array} 
\right. 
\end{equation} 
 As shown in Proposition~\ref{ffjkl} below, such a transformation 
 ${\tau}$ satisfies Condition~\eqref{djdkd1.0}. 
 The next figure shows an example of behaviour 
 such a transformation, with a 
 finite set of points for simplicity of illustration. 
\\ 
\vspace{0.4cm} 
\begin{center} 
\begin{picture}(230,130)(-60,-90)
\linethickness{0.2pt}
\put(-140,40){\circle{5}} 
\put(-140,-35){\circle{5}} 
\put(-140,-35){\line(0,1){75}} 
\put(-140,40){\line(2,1){50}} 
\put(-90,64.5){\circle{5}} 
\put(-90,64.5){\line(3,-1){90}} 
\put(0,35){\circle{5}} 
\put(0,35){\line(-1,-3){35}} 
\put(-35,-70){\circle{5}} 
\put(-35,-70){\line(-3,1){105}} 
\put(-60,0){\circle*{5}} 
\put(-65,-40){\circle*{5}} 
\put(-60,30){\circle*{5}} 
\put(-90,10){\circle*{5}} 
\put(-117,-22){\circle*{5}} 
\put(-77,23){\circle*{5}} 
\put(-86,-46){\circle*{5}} 
\put(-114,36){\circle*{5}} 

\put(10,0){$\stackrel{\tau}{\vector(3,0){90}}$} 

\put(120,40){\circle{5}} 
\put(120,-35){\circle{5}} 
\put(120,-35){\line(0,1){75}} 
\put(120,40){\line(2,1){50}} 
\put(170,64.5){\circle{5}} 
\put(170,64.5){\line(3,-1){90}} 
\put(260,35){\circle{5}} 
\put(260,35){\line(-1,-3){35}} 
\put(225,-70){\circle{5}} 
\put(225,-70){\line(-3,1){105}} 
\put(220,0){\circle*{5}} 
\put(195,40){\circle*{5}} 
\put(200,-30){\circle*{5}} 
\put(170,-10){\circle*{5}} 
\put(143,22){\circle*{5}} 
\put(183,23){\circle*{5}} 
\put(174,46){\circle*{5}} 
\put(146,-36){\circle*{5}} 
\end{picture}
\end{center} 
\vspace{-0.5cm} 
 
\begin{prop} 
\label{ffjkl} 
 The mapping ${\tau} : \Omega^X \times X \longrightarrow X$ 
given in \eqref{adfdsfg} satisfies Condition~\eqref{djdkd1.0}. 
\end{prop} 
\begin{Proof} 
 Let $x_1,\ldots ,x_m \in X$. Clearly, we can assume that some $x_i$ 
 lies outside of ${\cal C}(\omega)={\cal C}(\omega_e)$, otherwise 
 \begin{eqnarray*}
D_{x_i} \tau ( x_j , \omega)
&=&\tau(x_j,\omega\cup\{x_i\})-\tau(x_j,\omega)
=\tau(x_j,(\omega\cup\{x_i\})_e)-\tau(x_j,\omega_e)\\
&=&\tau(x_j,\omega_e)-\tau(x_j,\omega_e)= 0 
\end{eqnarray*}
 for all $i,j=1,\ldots ,m$. 
 Similarly, we can assume that ${\cal C}(\omega \cup \{ x_1,\ldots ,x_m \})$ 
 has at least one extremal point $x_i \in \{ x_1,\ldots ,x_m \}$. 
\\ 
 
 Now we have 
$$ 
 \tau ( x_i , \omega \cup \eta ) = 
 \tau ( x_i , \omega ) = x_i 
$$ 
 for all $\eta \subset \{ x_1,\ldots , x_m\}$, 
 hence 
$$ 
 D_\Theta \tau ( x_i , \omega ) = 0, 
$$ 
 for all $\Theta \subset \{ x_1,\ldots , x_m\}$, 
 due to the following consequence of \eqref{eq:Dmult}
\begin{eqnarray*} 
 D_\Theta 
 \tau ( x_i , \omega ) 
 & = & 
 \sum_{\eta \subset \Theta } 
 (-1)^{ | \Theta | + 1 - | \eta |} 
 \tau ( x_i , \omega \cup \eta ) 
\\ 
 & = & 
 \tau ( x_i , \omega ) 
 \sum_{\eta \subset \Theta } 
 (-1)^{ | \Theta | + 1 - | \eta |} 
\\ 
 & = & 
 \tau ( x_i , \omega \cup \eta ) 
 ( 1 - 1 )^{|\Theta |+1} 
\\ 
 & = & 0, 
\end{eqnarray*} 
 where the summation above holds over all (possibly empty) 
 subset $\eta$ of $\Theta$. 
As a consequence, one factor of \eqref{djdkd1.0} necessarily vanishes. 
\end{Proof} 
\section{Moment identities} 
\label{sec:moment}
From the previous factorial moment identities, we can 
recover some recently obtained 
moment identities for Poisson stochastic integrals 
with random integrands, 
cf. \cite{momentpoi}, 
and their extensions to point processes, 
cf. \cite{flint}. 
 Let 
\begin{equation}
\label{eq:Stirling}
 S(n,k) 
 = 
 \frac{1}{k!} 
 \sum_{d_1+\cdots + d_k=n} 
 \frac{n!}{d_1! \cdots d_k!} 
\end{equation}
 denote the Stirling number of the second kind, 
 i.e. the number of partitions of a set of $n$ 
 objects into $k$ non-empty subsets, 
 cf. also Relation~(3) page 2 of \cite{sloane}. 
 As a consequence we recover the following elementary moment identity
 from Proposition~\ref{klddd2}. 
\begin{lemma} 
\label{djkld} 
 Let $A=A(\omega)$ be a random set. 
 We have 
$$ 
 \Bbb{E} \left[ 
 F\ N(A)^n 
 \right] 
 = 
 \sum_{k=0}^n 
 S(n,k) 
 \Bbb{E} \left[ 
 \int_{X^k } 
 \varepsilon^+_{\eufrak{x}_k}( 
 F 
 \ind_{A} (x_1)  \cdots  \ind_{A} (x_k ) 
 ) 
\ \hat{\sigma} ( dx_1 , \ldots , dx_k ) 
 \right] 
. 
$$ 
\end{lemma} 
\begin{Proof} 
 This result is a direct consequence of 
 Proposition~\ref{klddd2} and the relation 
\begin{equation}
\label{eq:moment-factorial}
 \Bbb{E} [ X^n] = \sum_{k=1}^n S(n,k) \mu_k^f(X)
, 
\end{equation}
 between the moments and the factorial moments $\mu_k^f(X)$ 
 of a random variable $X$. 
This relation follows from the classical identity 
\begin{equation*} 
\label{djklddd} 
 x^n 
 = 
 \sum_{k=0}^n 
 S(n,k)\
 x (x-1) \cdots (x-k+1) 
,
\end{equation*} 
 cf. e.g. \cite{gerstenkorn} or page 72 of \cite{fristedt}. 
\end{Proof} 
 More generally, 
 Lemma~\ref{djkld} allows us to recover the 
 following moment identity, cf. 
 Theorem~3.1 of \cite{flint}, and 
 Proposition~3.1 of \cite{momentpoi} for the 
 Poisson case. 
\begin{prop} 
\label{djld} 
 Let $u : X \times \Omega^X \longrightarrow \real$ be 
 a (measurable) process. We have 
\begin{equation} 
\label{eq:momentu}
 \Bbb{E} \left[ 
 \left( \int_X u(x,\omega )\ \omega (dx) \right)^n 
 \right] 
 = 
 \sum_{k=1}^n 
 \sum_{B^n_1,\ldots , B^n_k } 
 \Bbb{E} \left[ 
 \int_{X^k} 
 \varepsilon^+_{\eufrak{x}_k} 
 \left( u(x_1,\cdot)^{|B^n_1|}\cdots u(x_k, \cdot)^{|B^n_k|} \right) 
 \ \hat{\sigma} ( d \eufrak{x}_k ) \right] 
\end{equation} 
 where the sum runs over the partitions 
 $B^n_1,\ldots , B^n_k$ of $\{ 1 , \ldots , n \}$, 
 for any $n\geq 1$ such that 
 all terms are integrable. 
\end{prop} 
\begin{Proof} 
 First we establish \eqref{eq:momentu} for 
 simple processes of the form 
 $u(x,\omega)= \sum_{i=1}^p F_i(\omega)\ind_{A_i(\omega)} (x)$ 
 with  a.s. disjoint random sets $A_i(\omega)$, $1\leq i\leq p$. 
 Applying Lemma~\ref{djkld} inductively we have 
\begin{eqnarray} 
\nonumber
\lefteqn{ 
 \Bbb{E} \left[ 
 \left( \sum_{i=1}^p F_i \int_X {\bf 1}_{A_i} (x)\ \omega (dx) \right)^n 
 \right] 
 = 
 \Bbb{E} \left[ 
 \left( \sum_{i=1}^p F_i N ( A_i ) \right)^n 
 \right] 
} 
\\ 
\nonumber
 & = & 
 \sum_{ n_1 + \cdots + n_p = n 
 \atop n_1,\ldots ,n_p \geq 0 
 } 
 \frac{n!}{n_1!\cdots n_p!} 
 \Bbb{E} \big[ 
 \left( 
 F_1 N ( A_1) 
 \right)^{n_1} 
 \cdots 
 \left( 
 F_p N ( A_p ) 
 \right)^{n_p} 
 \big] 
\\ 
\nonumber
 & = & 
 \sum_{ n_1 + \cdots + n_p = n 
 \atop n_1,\ldots ,n_p \geq 0 
 } 
 \frac{n!}{n_1!\cdots n_p!} 
 \sum_{k_1=0}^{n_1} 
 \cdots 
 \sum_{k_p=0}^{n_p} 
 S(n_1,k_1) 
 \cdots 
 S(n_p,k_p) 
\\
\nonumber
&&
 \Bbb{E} \left[ 
 \int_{X^{k_1+\cdots +k_p} }
 \varepsilon^+_{\eufrak{x}_{k_1+\cdots +k_p}} 
 \left( 
 F_1^{n_1} 
 \cdots 
 F_p^{n_p} 
 \ind_{A_1^{k_1}}\otimes \cdots \otimes
 \ind_{A_p^{k_p}} (x_1, \ldots, x_{k_1+\cdots +k_p} 
 ) 
 \right) 
 \right. 
\\ 
\nonumber
&& \hskip 5cm 
 \ \hat c(\{x_1, \ldots, x_{k_1+\cdots +k_p} \},\omega) \ 
 \sigma (dx_1) \cdots  \sigma (dx_{k_1+\cdots +k_p} ) \Big] \\
\nonumber
 & = & 
 \sum_{m=0}^n 
 \sum_{ n_1 + \cdots + n_p = n 
 \atop n_1,\ldots ,n_p \geq 0 
 } 
 \frac{n!}{n_1!\cdots n_p!} 
 \sum_{k_1+\cdots + k_p = m 
 \atop 1 \leq k_1 \leq n_1, \ldots , 1 \leq k_p \leq n_p }
 S(n_1,k_1) 
 \cdots 
 S(n_p,k_p) 
\\
\nonumber
&&
 \Bbb{E} \left[ 
 \int_{X^m} 
 \varepsilon^+_{\eufrak{x}_m} 
 \left( 
 F_1^{n_1} 
 \cdots 
 F_p^{n_p} 
 \ind_{A_1^{k_1}} \otimes \cdots \otimes
 \ind_{A_p^{k_p}} 
 (x_1, \ldots, x_m ) 
 \right) 
 \ \hat c(\{x_1, \ldots, x_m \},\omega) \ 
 \sigma (d\eufrak{x}_m ) \right] 
 \\
\nonumber 
 & = & 
 \sum_{m=0}^n 
 \sum_{ n_1 + \cdots + n_p = n 
 \atop n_1,\ldots ,n_p \geq 0 
 } 
 \frac{n!}{n_1!\cdots n_p!} 
 \sum_{I_1\cup \cdots\cup I_p=\{1, \cdots, m\} 
 \atop |I_1| \leq n_1, \ldots , |I_p| \leq n_p }
 S(n_1,|I_1|) 
 \cdots 
 S(n_p,|I_p|) 
 \frac{|I_1|!\cdots|I_p|!}{m!}
\\
\nonumber 
&&
 \Bbb{E} \left[ 
 \int_{X^m} 
 \varepsilon^+_{\eufrak{x}_m} 
 \left( 
 F_1^{n_1} 
 \cdots 
 F_p^{n_p} 
 \prod_{j\in I_1} \ind_{A_1}(x_j) \cdots \prod_{j\in I_p}  \ind_{A_p}(x_j) 
 \right) 
 \ \hat c(\{x_1, \ldots, x_m \},\omega) \ 
 \sigma (d\eufrak{x}_m ) \right] 
\\ 
\label{eq:change-var}
\\ 
\nonumber 
 & = & 
 \sum_{m=0}^n 
 \sum_{ P_1\cup \cdots\cup P_m = \{ 1 , \ldots , n \} } 
 \sum_{ i_1,\ldots ,i_m = 1}^p 
 \Bbb{E} \left[ 
 \int_{X^m} 
 \varepsilon^+_{ \eufrak{x}_m} 
 \left( F_{i_1}^{|P_1|} 
 {\bf 1}_{A_{i_1}} (x_1) 
 \cdots 
 F_{i_m}^{|P_m|} 
 {\bf 1}_{A_{i_m}} (x_m) \right)
  \ \hat{\sigma} ( d\eufrak{x}_m ) 
 \right], 
\\ 
\label{eq:Stir1a} 
\end{eqnarray}
 where in \eqref{eq:change-var} we made changes of variables in the integral 
 and, in \eqref{eq:Stir1a}, we used the combinatorial identity of Lemma \ref{lemme:Stir1} below 
with $\alpha_{i,j}=\ind_{A_i}(x_j)$, $1\leq i\leq p, 1\leq j\leq m$, and $\beta_i=F_i$. 
 The proof is concluded by 
 using the disjunction of the $A_i$'s in \eqref{eq:Stir2}, 
 as follows:
\begin{eqnarray} 
\nonumber
\lefteqn{ 
 \Bbb{E} \left[ 
 \left( \sum_{i=1}^p F_i \int_X {\bf 1}_{A_i} (x)\ \omega (dx) \right)^n 
 \right] 
} 
\\ 
\nonumber
& = & 
 \sum_{m=0}^n 
 \sum_{ P_1\cup \cdots\cup P_m = \{ 1 , \ldots , n \} } 
 \Bbb{E} \left[ 
 \int_{X^m} 
 \varepsilon^+_{ \eufrak{x}_m} \left(
 \sum_{i=1}^p \left( F_i^{|P_1|} 
 {\bf 1}_{A_i} (x_1) 
 \right)
 \cdots 
 \sum_{i=1}^p 
 \left( 
 F_i^{|P_m|} {\bf 1}_{A_i} (x_m) 
 \right) 
\right)
 \ \hat{\sigma} ( d\eufrak{x}_m ) 
 \right] 
\\
\nonumber
& = & 
 \sum_{m=0}^n 
 \sum_{ P_1\cup \cdots\cup P_m = \{ 1 , \ldots , n \} } 
 \Bbb{E} \left[ 
 \int_{X^m} 
 \varepsilon^+_{ \eufrak{x}_m} \left(
 \left( 
 \sum_{i=1}^p F_i {\bf 1}_{A_i} (x_1) 
 \right)^{|P_1|} 
 \cdots 
 \left( 
 \sum_{i=1}^p F_i {\bf 1}_{A_i} (x_m) 
 \right)^{|P_m|} \right)
 \ \hat{\sigma} ( d\eufrak{x}_m ) 
 \right].\\
\label{eq:Stir2}
\end{eqnarray} 
 The general case is obtained by approximating $u(x,\omega )$ 
 with simple processes. 
\end{Proof} 
Using \eqref{eq:espilon-mult}, we can also write 
$$ 
 \Bbb{E} \left[ 
 \left( \int_X u(x,\omega )\ \omega (dx) \right)^n 
\right] 
 = 
 \sum_{k=1}^n 
 \sum_{B^n_1,\ldots , B^n_k } 
 \sum_{\Theta \subset \{1,\ldots ,k\} } 
 \Bbb{E} \left[ 
 \int_{X^k} 
 D_{\Theta} 
 ( 
 u^{|B^n_1|}_{x_1} \cdots 
 u^{|B^n_k|}_{x_k} ) \ \hat{\sigma} ( d \eufrak{x}_k ) \right] 
. 
$$ 
 The next lemma has been used above in the proof of 
 Proposition~\ref{djld}. 
\begin{lemma} 
\label{lemme:Stir1} 
Let $m,n,p\in\Bbb{N}$, $(\alpha_{i,j})_{1\leq i\leq p, 1\leq j\leq m}$
and $\beta_1 , \ldots , \beta_p \in\real$. 
We have 
\begin{eqnarray} 
\nonumber 
&&\sum_{ n_1 + \cdots + n_p = n 
 \atop n_1,\ldots ,n_p \geq 0 
 } 
 \frac{n!}{n_1!\cdots n_p!} 
 \sum_{I_1\cup\cdots\cup I_p =\{1,\ldots, m\} 
 \atop |I_1|\leq n_1, \ldots , |I_p| \leq n_p }
 S(n_1,|I_1|) 
 \cdots 
 S(n_p,|I_p|) \times
 \\
 \nonumber
 && \hskip 5cm
 \frac{|I_1|!\cdots |I_p|!}{m!}
  \beta_1^{n_1} 
\Big(\prod_{j\in I_1}\alpha_{1,j}\Big) 
 \cdots 
 \beta_p^{n_p} 
 \Big(\prod_{j\in I_p}\alpha_{p,j}\Big) 
\\ 
\label{ddd} 
 & = & 
 \sum_{ P_1\cup \cdots\cup P_m = \{ 1 , \ldots , n \} } 
 \sum_{ i_1,\ldots ,i_m = 1}^p 
 \beta_{i_1}^{|P_1|} 
 \alpha_{i_1,1} 
 \cdots 
 \beta_{i_m}^{|P_m|} 
 \alpha_{i_m,m}.
\end{eqnarray} 
\end{lemma} 
\begin{Proof} 
Observe that \eqref{eq:Stirling} ensures
\begin{eqnarray*}
 S(n,|I|) 
 \beta^{n} 
\Big(\prod_{j\in I}\alpha_j\Big) 
&=&\sum_{ 
 \bigcup_{a\in I} P_a = \{ 1 , \ldots , n \} } 
 \prod_{j\in I} \big(\alpha_j\beta^{|P_j|}\big) 
\end{eqnarray*}
for all $\alpha_j$, $j\in I$, $\beta\in\real$, $n\in \Bbb{N}$. 
 We have 
\begin{eqnarray*} 
\lefteqn{ 
 \sum_{ n_1 + \cdots + n_p = n 
 \atop n_1,\ldots ,n_p \geq 0 
 } 
 \frac{n!}{n_1!\cdots n_p!} 
 \sum_{I_1\cup\cdots\cup I_p=\{1, \ldots, m\} 
 \atop |I_1| \leq n_1, \ldots , |I_p| \leq n_p }
 S(n_1,|I_1|) 
 \cdots 
 S(n_p,[I_p|)
 }
  \\
&&\hskip 5cm  \frac{|I_1|!\cdots |I_p|!}{m!}
 \beta_1^{n_1} 
 \big(\prod_{j\in I_1}\alpha_{1,j}\Big)
 \cdots 
 \beta_p^{n_p} 
\Big(\prod_{j\in I_p}\alpha_{p,j}\Big)  
\\
&=& \sum_{ n_1 + \cdots + n_p = n 
 \atop n_1,\ldots ,n_p \geq 0 
 } 
 \frac{n!}{n_1!\cdots n_p!}\sum_{I_1\cup \cdots \cup I_p =\{1, \ldots, m\} 
 \atop |I_1| \leq n_1, \ldots , |I_p| \leq n_p }\frac{|I_1|!\cdots |I_p|!}{m!}
\\
&&
\Big(\sum_{ ~~\bigcup_{a\in I_1} P_a^1 = \{ 1 , \ldots , n_1 \} } 
 \prod_{j_1\in I_1} \big(\alpha_{1,j_1}\beta_1^{|P_{j_1}^1|}\big) \Big)
\cdots
\Big(\sum_{ ~~\bigcup_{a\in I_p} P_a^p = \{ 1 , \ldots , n_p \} } 
 \prod_{j_p\in J_p} \big(\alpha_{p,j_p}\beta_p^{|P_{j_p}^p|}\big) \Big)
\\
&=& \sum_{ n_1 + \cdots + n_p = n 
 \atop n_1,\ldots ,n_p \geq 0 
 } 
 \frac{n!}{n_1!\cdots n_p!}  
\sum_{I_1\cup \cdots \cup I_p = \{1, \ldots, m\} 
 \atop |I_1| \leq n_1, \ldots , |I_p| \leq n_p }
\sum_{ ~~\bigcup_{a\in I_1} P_a^1 = \{ 1 , \ldots , n_1 \} }
\cdots \sum_{ ~~\bigcup_{a\in I_p} P_a^p = \{ 1 , \ldots , n_p \} }
\\
&&\hskip 5cm \frac{|I_1|!\cdots |I_p|!}{m!}
\prod_{l=1}^p \prod_{j_l\in I_l} \big(\alpha_{l,j_l}\beta_l^{|P_{j_l}^l|}\big)
\\
 & = & 
 \sum_{ n_1 + \cdots + n_p = n 
 \atop n_1,\ldots ,n_p \geq 0 
 } 
 \frac{n!}{n_1!\cdots n_p!}  
\sum_{I_1\cup \cdots \cup I_p = \{1, \ldots, m\} 
 \atop |I_1| \leq n_1, \ldots , |I_p| \leq n_p }
\sum_{ ~~\bigcup_{a\in I_1} P_a^1 = \{ 1 , \ldots , n_1 \} }
\cdots \sum_{ ~~\bigcup_{a\in I_p} P_a^p = \{ 1 , \ldots , n_p \} }
\\
&&\hskip 5cm \frac{|I_1|!\cdots |I_p|!}{m!}
\prod_{l=1}^p \prod_{j_l\in I_l} \alpha_{l,j_l} 
\prod_{l=1}^p \prod_{j_l\in I_l} \beta_l^{|P_{j_l}^l|} 
\\
 & = & 
 \sum_{ n_1 + \cdots + n_p = n 
 \atop n_1,\ldots ,n_p \geq 0 
 } 
 \frac{n!}{n_1!\cdots n_p!}  
 \sum_{k_1 + \cdots + k_p = m 
 \atop 1 \leq k_1 \leq n_1, \ldots , 1 \leq k_p \leq n_p }
 \sum_{ i_1,\ldots ,i_m = 1}^p 
\\
&& 
 \sum_{ P_1^1\cup \cdots\cup P_{k_1}^1 = \{ 1 , \ldots , n_1 \} }
 \cdots 
 \sum_{ P_{k_1+\cdots +k_{p-1}+1}^p\cup \cdots\cup P_{k_1+\cdots + k_p}^p = \{ 1 , \ldots , n_p \} }
 \prod_{j=1}^m \big( \alpha_{i_j,j} \beta_{i_j}^{|P_j^{i_1}|+\cdots + | P_j^{i_m} | } \big) 
\\
 & = & 
 \sum_{ P_1\cup \cdots\cup P_m = \{ 1 , \ldots , n \} } 
 \sum_{ i_1,\ldots ,i_m = 1}^p 
 \beta_{i_1}^{|P_1|} 
 \alpha_{i_1,1} 
 \cdots 
 \beta_{i_m}^{|P_m|} 
 \alpha_{i_m,m}  
, 
\end{eqnarray*} 
 by a reindexing 
 of the summations and the fact that the reunions of the partitions 
$P_1^j,\ldots,P_{|I_j|}^j$, $1\leq j\leq p$, of disjoint $p$ subsets of $\{1,\ldots, m\}$ run the partition of $\{1, \ldots, m\}$ when we take into account the choice of the $p$ subsets and the possible length $k_j$, $1\leq j\leq p$, of the partitions. 
\end{Proof} 
 Note that the combinatorial result of Lemma~\ref{lemme:Stir1} 
 can also be shown in a probabilistic way 
 when $\alpha_{i,j} = \alpha_i$, $1\leq i\leq p$, $1\leq j\leq m$. 
 Recall the relation \eqref{eq:moment-factorial} between standard 
 moments and factorial moments. 
 From \eqref{eq:2} we have 
\begin{eqnarray} 
\nonumber 
\lefteqn{ 
 \sum_{m=0}^n 
 \lambda^m 
 \sum_{ n_1 + \cdots + n_p = n 
 \atop n_1,\ldots ,n_p \geq 0 
 } 
 \frac{n!}{n_1!\cdots n_p!} 
 \sum_{k_1+\cdots + k_p = m 
 \atop k_1 \leq n_1, \ldots , k_p \leq n_p }
 S(n_1,k_1) 
 \cdots 
 S(n_p,k_p) 
 \beta_1^{n_1} 
 \alpha_1^{k_1} 
 \cdots 
 \beta_p^{n_p} 
 \alpha_p^{k_p} 
} 
\\ 
\nonumber 
 & = & 
 \sum_{ n_1 + \cdots + n_p = n 
 \atop n_1,\ldots ,n_p \geq 0 
 } 
 \frac{n!}{n_1!\cdots n_p!} 
 \sum_{k_1=0}^{n_1} 
 S(n_1,k_1) 
 ( \lambda \alpha_1 )^{k_1} 
 \cdots 
 \sum_{k_p=0}^{n_p} 
 S(n_p,k_p) 
 ( \lambda \alpha_p )^{k_p} 
 \beta_1^{n_1} 
 \cdots 
 \beta_p^{n_p} 
\\ 
\nonumber 
 & = & 
 \sum_{ n_1 + \cdots + n_p = n 
 \atop n_1,\ldots ,n_p \geq 0 
 } 
 \frac{n!}{n_1!\cdots n_p!} 
 \beta_1^{n_1} 
 \cdots 
 \beta_p^{n_p} 
 \Bbb{E}[ 
 Z_{\lambda \alpha_1}^{n_1} 
 \cdots 
 Z_{\lambda \alpha_p}^{n_p} 
 ] 
\\ 
\nonumber 
 & = & 
 \Bbb{E}\left[ 
 \left( 
 \sum_{i=1}^p  \beta_i Z_{\lambda \alpha_i} 
 \right)^n 
 \right] 
\\ 
\label{3} 
 & = & 
 \sum_{m=0}^n 
 \lambda^m 
 \sum_{ P_1\cup \cdots\cup P_m = \{ 1 , \ldots , n \} } 
 \sum_{ i_1,\ldots ,i_m = 1}^p 
 \beta_{i_1}^{|P_1|} 
 \alpha_{i_1} 
 \cdots 
 \beta_{i_m}^{|P_m|} 
 \alpha_{i_m}, 
\end{eqnarray} 
 since by \eqref{eq:2} the moment of order $n_i$ of $Z_{\lambda \alpha_i}$ 
 is given by 
$$ 
 \Bbb{E} \left[ Z_{\lambda \alpha_i}^{n_i} \right] 
 = \sum_{k=0}^{n_i} S(n_i,k) (\lambda \alpha_i)^k. 
$$ 
 The above relation \eqref{3} 
 being true for all $\lambda$, this implies 
 \eqref{ddd} for this choice of $\alpha_{i,j}$'s. 
\section{Poisson case} 
\label{sec:invariance}
 In the Poisson case, we have $c(x,\omega ) =1$ 
 and the results of the previous sections 
 specialize immediately 
 to new factorial moment identities for 
 Poisson point processes with intensity $\sigma (dx)$. 
 For any random set $A=A (\omega )$ 
 and sufficiently integrable random variable $F$, we have 
\begin{eqnarray*} 
 \Bbb{E} \left[ 
 F\ N(A)_{(n)} 
 \right] 
 & = & 
 \Bbb{E} \left[ 
 \int_{X^n} 
\varepsilon^+_{\eufrak{x}_n}( 
 F 
 \ind_{A} (x_1)  \cdots  \ind_{A} (x_n)) 
\ \sigma (dx_1) \cdots  \sigma (dx_n) \right] 
, 
\end{eqnarray*}  
 $n\geq 1$. 
 For all almost surely disjoint random sets $A_i(\omega)$, $1\leq i\leq p$, 
 and sufficiently integrable random variable $F$, we have 
\begin{eqnarray*}  
\lefteqn{ 
 \Bbb{E} \left[ 
 F\ N(A_1)_{(n_1)} 
 \cdots N(A_p)_{(n_p)} 
 \right] 
} 
\\ 
\nonumber 
 & = & 
 \Bbb{E} \left[ 
 \int_{X^n} 
 \varepsilon^+_{\eufrak{x}_n}( 
 F 
(\ind_{A_1^{n_1}}\otimes \cdots \otimes
 \ind_{A_p^{n_p}})(x_1, \ldots, x_n) 
 ) 
\ \sigma (dx_1) \cdots  \sigma (dx_n) \right] 
, 
\end{eqnarray*}
 with $n=n_1+\cdots+n_p$. 
 In addition, we have 
 the following proposition whose proof is similar to that 
 of Proposition~\ref{djhjklddd111} 
 although it cannot be obtained as a direct 
 consequence of Proposition~\ref{djhjklddd111} 
 and it cannot be stated in the (non-Poisson) 
 point process setting. 
\begin{prop} 
\label{djklddd111} 
 Consider $A_1 (\omega ), \ldots, A_p (\omega )$ 
 a.s. disjoint random sets such that 
 $\sigma (A_i (\omega ))$ is deterministic, 
 $i=1,\ldots ,p$, and 
\begin{equation} 
\label{djdkd1.1}
 D_{\Theta_1} 
 {\bf 1}_{A_i ( \omega )}(x_1) 
 \cdots 
 D_{\Theta_m} 
 {\bf 1}_{A_i ( \omega )}(x_m) 
 = 
 0, 
\end{equation} 
 for every family $\{ \Theta_1, \ldots , \Theta_m \}$ 
 of (non empty) subsets such that 
 $\Theta_1 \cup \cdots \cup \Theta_m = \{1,\ldots ,m \}$, 
 all $x_1,\ldots ,x_m \in X$, all $\omega\in\Omega^X$. 
 Then the family 
$$ 
\big( N ( A_1 ) , \ldots , N(A_p )\big) 
$$ 
 is a vector of independent Poisson random 
 variables with parameters 
 $\sigma (A_1 ), \ldots , \sigma (A_p )$. 
\end{prop} 
\begin{Proof} 
 Let $n=n_1+\cdots+n_p$. 
Under Condition~\eqref{djdkd1.1}, 
Lemma~\ref{l12.32} and \eqref{eq:DDD} show that 
\begin{equation}
\label{eq:cons2lemma}
 D_{x_1} \cdots D_{x_k} \big( 
 {\bf 1}_{A_{i_1} ( \omega )}(x_1) 
 \cdots 
 {\bf 1}_{A_{i_k} ( \omega )}(x_k) 
\big) = 0, 
 \qquad 
 x_1, \ldots , x_k \in X, 
\end{equation}
 for all $i_1,\ldots ,i_k \in \{ 1,\ldots ,n \}$ and 
 $\omega \in \Omega$. 
 Since in addition $\sigma ( A_i )$ is deterministic, 
 $i=1,\ldots ,p$, then by 
 \eqref{eq:DThetajoint} with $F=1$ we obtain 
\begin{eqnarray} 
\nonumber
\lefteqn{ 
\! \! \! \! \! \! \! \! \! \! \! \! \! \! \! \! \! \! \! \! \! \! \! \! \! \! \! \! \! \! 
 \Bbb{E} \left[ 
 N(A_1)_{(n_1)} 
 \cdots N(A_p)_{(n_p)} 
 \right] 
 = 
\sum_{\Theta \subset \{1,\ldots ,n\} } 
\Bbb{E} \left[ 
 \int_{X^n} 
 D_{\Theta} 
 \big( 
 ( 
 \ind_{A_1^{n_1}}\otimes \cdots \otimes
 \ind_{A_p^{n_p}} ) (\eufrak{x}_n ) 
 \big)\ 
 {\sigma}^n ( d\eufrak{x}_n ) 
 \right]
} 
\\
\nonumber
&=& 
\Bbb{E} \left[ 
 \int_{X^n} 
 D_{x_1} \cdots D_{x_n} 
 \big( 
 ( 
 \ind_{A_1^{n_1}}\otimes \cdots \otimes
 \ind_{A_p^{n_p}})(\eufrak{x}_n ) 
\big)\   
 {\sigma}^n ( d\eufrak{x}_n ) 
 \right]
\\ 
\nonumber
 & & + 
\sum_{\Theta \subsetneq \{1,\ldots ,n\} } 
\Bbb{E} \left[ 
 \int_{X^n} 
 D_{\Theta} 
 \big( 
 (\ind_{A_1^{n_1}}\otimes \cdots \otimes
 \ind_{A_p^{n_p}})(\eufrak{x}_n ) 
\big)\ 
 {\sigma}^n ( d\eufrak{x}_n ) 
 \right]
\\ 
\label{eq:generic}
 & = & 
\sum_{\Theta \subsetneq \{1,\ldots ,n\} } 
\Bbb{E} \left[ 
 \int_{X^n} 
 D_{\Theta} 
 \big( 
 (\ind_{A_1^{n_1}}\otimes \cdots \otimes 
 \ind_{A_p^{n_p}})(\eufrak{x}_n ) 
\big)\ 
 {\sigma}^n ( d\eufrak{x}_n ) 
 \right]
\end{eqnarray} 
using \eqref{eq:cons2lemma}. 
 Next, without loss of generality the generic, term of \eqref{eq:generic} can be reduced to the term with  
 $\Theta = \{ 1 , \ldots , n-1 \}$ and using \eqref{eq:DDD}, we have 
 \begin{eqnarray} 
\nonumber
\lefteqn{ 
 \Bbb{E} \left[ 
 \int_{X^n} 
 D_{x_1} \cdots D_{x_{n-1}} 
\big( 
 (\ind_{A_1^{n_1}}\otimes \cdots \otimes
 \ind_{A_p^{n_p}})(\eufrak{x}_n ) 
\big)\ 
 {\sigma}^n ( d\eufrak{x}_n ) 
 \right] 
} 
\\ 
\nonumber
&=&
\Bbb{E} \left[ 
 \int_{X^n} 
 \sum_{\Theta_1\cup\cdots\cup \Theta_n = \{1, \ldots, n-1\}} 
 \prod_{k=1}^p\prod_{j=1}^{n_k} D_{\Theta_{n_1+\cdots+n_{k-1}+j}}\ind_{A_k}(x_{n_1+\cdots+n_{k-1}+j})
 \ 
 {\sigma}^n ( d\eufrak{x}_n ) 
 \right] \\
\label{eq:sigmaAp}
&=&
\Bbb{E} \left[ 
 \int_{X^{n-1}} 
 \sum_{\Theta_1\cup\cdots\cup \Theta_n = \{1, \ldots, n-1\}} 
\prod_{k=1}^{p-1}\prod_{j=1}^{n_k} D_{\Theta_{n_1+\cdots+n_{k-1}+j}}
\ind_{A_k}(x_{n_1+\cdots+n_{k-1}+j}) \right .\\
\nonumber
&&\left .
\times \prod_{j=1}^{n_p-1} D_{\Theta_{n_1+\cdots+n_{p-1}+j}}\ind_{A_p}(x_{n_1+\cdots+n_{p-1}+j})  
\int_X D_{\Theta_n}\ind_{A_p}(x_n)\ {\sigma} ( dx_n )\ {\sigma}^{n-1} ( d\eufrak{x}_{n-1} ) 
 \right] 
 \\
\nonumber
 &=&
\Bbb{E} \left[ 
 \int_{X^{n-1}} 
 \sum_{\Theta_1\cup\cdots\cup \Theta_{n-1} = \{1, \ldots, n-1\}} 
\prod_{k=1}^{p-1}\prod_{j=1}^{n_k} D_{\Theta_{n_1+\cdots+n_{k-1}+j}}
\ind_{A_k}(x_{n_1+\cdots+n_{k-1}+j}) \right .\\
\label{eq:tech36}
&&\left .
\times \prod_{j=1}^{n_p-1} D_{\Theta_{n_1+\cdots+n_{p-1}+j}}\ind_{A_p}(x_{n_1+\cdots+n_{p-1}+j})  
\sigma(A_p)\ {\sigma}^{n-1} ( d\eufrak{x}_{n-1} ) 
 \right] 
\\
\nonumber 
  & = & 
 \sigma( A_p) 
 \Bbb{E} \left[ 
 \int_{X^{n-1}} 
 D_{x_1} \cdots D_{x_{n-1}} 
 ( 
 (\ind_{A_1^{n_1}}\otimes \cdots \otimes
 \ind_{A_p^{n_p-1}})(\eufrak{x}_{n-1} ) 
 ) 
\ {\sigma}^{n-1} ( d\eufrak{x}_{n-1} ) 
 \right] 
\end{eqnarray} 
where \eqref{eq:tech36} comes from the fact that in \eqref{eq:sigmaAp} only the term with $\Theta_n=\emptyset$ is not zero
since 
$$
\int_X D_\Theta\ind_{A_p}(x)\ \sigma(dx)
=D_\Theta\left(\int_X\ind_{A_p}(x)\sigma(dx)\right)
=D_\Theta\big(\sigma(A_p)\big)=0
$$
using $\sigma ( A_p )$ is deterministic. 
 Finally, by applying the above argument recursively we obtain  
\begin{eqnarray*} 
 \Bbb{E} \left[ 
 N(A_1)_{(n_1)} 
 \cdots N(A_p)_{(n_p)} 
 \right] 
 & = & 
 \Bbb{E} \left[ 
 \int_{X^n} 
 (\ind_{A_1^{n_1}}\otimes \cdots \otimes
 \ind_{A_p^{n_p}})(\eufrak{x}_n ) 
 {\sigma}^n ( d\eufrak{x}_n ) 
 \right] 
\\ 
 & = & 
 \sigma ( A_1 )^{n_1} \cdots  \sigma ( A_p )^{n_p}, 
\end{eqnarray*} 
 $n_1,\ldots ,n_p \geq 1$, which characterizes the Poisson distribution 
 with parameters 
$$ 
 (\sigma(A_1), \ldots, \sigma(A_p)). 
$$ 
\end{Proof} 
 As a consequence, we recover the following 
 invariance result for Poisson measures when $(X,\sigma ) = (Y,\mu )$,
 where Condition~\eqref{cyclic3} below is interpreted as in 
 \eqref{eq:interpretation} above. 
\begin{theorem} 
\label{c10} 
 Let $\tau : \Omega^X \times X \to Y$ 
 be a random transformation such that $\tau ( \cdot , \omega ): X \to Y$ 
 maps $\sigma$ to $\mu$ for all $\omega \in \Omega^X$, 
 i.e. 
$$ 
 \sigma\circ \tau( \cdot , \omega )^{-1}  = \mu, 
 \qquad 
 \omega \in \Omega^X, 
$$ 
 and satisfying the condition 
\begin{equation} 
\label{cyclic3} 
 D_{\Theta_1} 
 \tau ( x_1 , \omega ) 
 \cdots 
 D_{\Theta_m} 
 \tau ( x_m , \omega ) 
 = 
 0, 
\end{equation} 
 for every family $\{ \Theta_1, \ldots , \Theta_m \}$ 
 of (non empty) subsets such that 
 $\Theta_1 \cup \cdots \cup \Theta_m = \{1,\ldots ,m \}$, 
 all $x_1,\ldots ,x_m \in X$, all $\omega\in\Omega^X$, 
 and all $i=1,\ldots ,p$. 
 Then $\tau_* : \Omega^X \to \Omega^Y$ defined by 
$$ 
 \tau_* \omega = \sum_{x\in \omega } \epsilon_{\tau ( x,\omega )}= \omega\circ \tau(\cdot, \omega)^{-1}, 
 \qquad 
 \omega \in \Omega^X, 
$$  
 maps $\pi_\sigma$ to $\pi_\mu$, i.e. 
$ 
 \tau_* \pi_\sigma 
$ 
 is the Poisson measure $\pi_\mu$ with intensity $\mu (dy)$ on $Y$. 
\end{theorem} 
\begin{Proof} 
 For any family $B_1,\ldots ,B_p$ of disjoint measurable subsets of 
 $Y$ with finite measure, we let 
 $A_i ( \omega ) = \tau^{-1} ( B_i , \omega ) \subset X$, i.e. 
 ${\bf 1}_{A_i} ( \cdot ) = {\bf 1}_{B_i} \circ \tau ( \cdot , \omega )$, 
 $i=1,\ldots ,p$, and by 
 Proposition~\ref{djklddd111}, we find that 
$$ 
 \omega \longmapsto 
 ( 
 \tau_* \omega ( B_1 ) , 
 \ldots , 
 \tau_* \omega ( B_p ) 
 ) 
 = 
 ( 
 \omega ( A_1) 
 ,\ldots, 
 \omega ( A_p ) 
 ) 
$$ 
 is a vector of independent Poisson random 
 variables with parameters $\mu (A_1 ), \ldots , \mu (A_p )$ 
 since $\sigma ( A_i ( \omega ) ) = 
 \sigma ( \tau^{-1} ( B_i , \omega ) ) 
 = \mu ( B_i )$ is deterministic, $i=1,\ldots ,p$, and 
\eqref{djdkd1.1} comes from the following consequence of \eqref{cyclic3}:
\begin{eqnarray*} 
\lefteqn{ 
 D_{\Theta_1} {\bf 1}_{A_{i_1} (\omega )} ( x_1 ) 
 \cdots 
 D_{\Theta_m} 
 {\bf 1}_{A_{i_m} (\omega )} ( x_m ) 
} 
\\ 
 & = & 
 D_{\Theta_1} 
 {\bf 1}_{B_{i_1}} ( \tau ( x_1 , \omega ) ) 
 \cdots 
 D_{\Theta_m} 
 {\bf 1}_{B_{i_m}} ( \tau ( x_m , \omega ) ) 
\\ 
 & = & 0. 
\end{eqnarray*} 
\end{Proof} 
 The example of random transformation given page \pageref{adfdsfg}
 at the end of Section~\ref{3.1} also satisfies 
 Condition~\eqref{cyclic3} in Theorem~\ref{c10}. 

\footnotesize 

\def\cprime{$'$} \def\polhk#1{\setbox0=\hbox{#1}{\ooalign{\hidewidth
  \lower1.5ex\hbox{`}\hidewidth\crcr\unhbox0}}}
  \def\polhk#1{\setbox0=\hbox{#1}{\ooalign{\hidewidth
  \lower1.5ex\hbox{`}\hidewidth\crcr\unhbox0}}} \def\cprime{$'$}

\end{document}